\documentclass[a4paper, 12pt, leqno]{amsart}
\usepackage{mathrsfs}
\usepackage{amsmath}
\usepackage{amscd}
\usepackage{amssymb}
\usepackage{amsthm}
\usepackage{pb-diagram}
\usepackage[all]{xy}

\numberwithin{equation}{section}

\def\Gg{{\mathfrak{g}}}

\def\BC{{\mathbb{C}}}

\def\BZ{{\mathbb{Z}}}
\def\CA{{\mathcal A}}
\def\CB{{\mathcal B}}

\def\CD{{\mathcal D}}

\def\CN{{\mathcal N}}

\def\vp{\varphi_0}

\title[The based ring of the lowest two-sided cell]{The based ring of the lowest two-sided cell of an affine Weyl group, III}

\author{Nanhua XI}
\address{Institute of Mathematics\\
Chinese Academy of Sciences\\
Beijing, 100190\\
and\\
School of Mathematical Sciences\\
 University of Chinese Academy of Sciences\\ Beijing 100049, China}
\email{nanhua@math.ac.cn}
\thanks{The work was partially supported by Natural Sciences Foundation
of China (No. 11321101).}

\begin{document}
\begin{abstract}
We show that Lusztig's homomorphism from an affine Hecke algebra to the direct summand of its asymptotic Hecke algebra corresponding to the lowest two-sided cell is related to the
 homomorphism constructed by Chriss and Ginzburg using equivariant K-theory by a matrix over the representation ring of the associated algebraic group.
\end{abstract}
\maketitle \setcounter{section}{0}
\label{sec:Intro} Let $W$ be the extended affine Weyl group associated to a simply connected simple algebraic group $G$ over $\BC$. Let $J_0$ be the based ring of the lowest two-sided cell of $W$ and let $H$ be the Hecke algebra of $W$ over the ring
$\CA=\BZ[v,v^{-1}]$ of  Laurent polynomials in an indeterminate $v$ with
integer coefficients. Lusztig defined an $\CA$-algebra homomorphism $\vp:H\to J_0\otimes_{\BZ}\CA$ (see [L3]). In [CG,  Corollary 5.4.34], Chriss and Ginzburg defined an $\CA$-algebra homomorphism $\psi_0:H\to J_0\otimes_{\BZ}\CA$ using equivariant K-theory. We will show that $\vp$ is essentially the conjugacy of $\psi_0$ by an invertible  $W_0\times W_0$ matrix with entries in the representation $R_{G\times\BC^*}$ of $G\times\BC^*$, see Theorem 2.5.

\section{Preliminaries}

\noindent{\bf 1.1.} Let $G$ be a simply connected simple algebraic
group over the complex number field $\BC$. The Weyl group $W_0$ of $G$ acts
naturally on the character group $X$ of a maximal tours of $G$. The
semidirect product $W=W_0\ltimes X$ with respect to the action is
called an (extended) affine Weyl group. Let $H$ be the associated
Hecke algebra over the ring $\CA=\BZ[v,v^{-1}]$ ($v$ an
indeterminate) with parameter $v^2$. Thus $H$ has an $\CA$-basis
$\{T_w\ |\ w\in W\}$ and its multiplication is defined by the
relations $(T_s-v^2)(T_s+1)=0$ if $s$ is a simple reflection and
$T_wT_u=T_{wu}$ if $l(wu)=l(w)+l(u)$, here $l$ is the length
function of $W$.
\def\gb{\mathfrak{b}}
\def\a{\alpha}

\smallskip

We shall write the operation in $X$ multiplicatively. Let $X^+=\{x\in X\,|\, l(w_0x)=l(w_0)+l(x)\}$ be the set of dominant weights in $X$, where $w_0$ is the longest element of $W_0$. For any $z\in X$, one can find dominant weights $x$ and $y$ such that $z=xy^{-1}$. Then define $\theta_z=v^{l(y)-l(x)}T_xT_y^{-1}$. It is known that $\theta_z$ is well defined, that is, it only depends on $z$ and independent of the choice of $x,y$. Moreover, the elements $\theta_zT_w,\ z\in X,\ w\in W_0$, form an $\CA$-basis of $H$ and the  the elements $T_w\theta_z,\  w\in W_0$, $z\in X,$ form an $\CA$-basis of $H$ as well. See [L1] or [L4].

\smallskip

The $\CA$-subalgebra $\Theta$ of $H$ generated by all  $\theta_x,\ x\in X$, is commutative and is isomorphic to the group algebra $\CA[X]$ of $X$. Let $\theta:\CA[X]\to \Theta$, $\sum_{x\in X}a_xx\to \sum_{x\in X}a_x\theta_x$ be the isomorphism, where $a_x\in\CA$ are 0 except for finitely many $x$ in $X$. For $\gamma\in\CA[X]$, we shall also write $\theta_{\gamma}$ for $\theta(\gamma)$.

Let $R$ be the root system. For simple reflection $s$ in $W_0$, we shall denote by $\alpha_s$ for the simple root which defines $s$. For $x\in X$,  we have (see [L5, 1.19])

\vskip.2cm

(a) $\displaystyle{T_s\theta_{x}=\theta_{s(x)}T_s-(v^2-1)\theta_{\frac{s(x)\a_s-x\a_s}{\a_s-1}}}$.

\vskip.3cm \noindent{\bf 1.2.} For $y,w$ in $W$, let $P_{y,w}$ be the Kazhdan-Lusztig polynomial. Set $C_w=v^{l(w)}\sum_{y\le w}(-1)^{l(w)-l(y)}P_{y,w}(v^{-2})T_y,$  where $\le$ is the Bruhat order on $W$. Then the elements $C_w,\ w\in W$, form an $\CA$-basis of $H$ and is called a Kazhdan-Lusztig basis of $H$.

Recall that $w_0$ is the longest element in $W_0$. We shall simply write $C$ for $C_{w_0}$. We have
$T_sC=-C$ for all simple reflections $s$ in $W_0$. Therefore, the $\CA$-linear map $\pi$: $\CA[X]\to HC$ defined by $x\to \theta_xC$ for all $x\in X$ is an isomorphism of $\CA$-module. Using the formula 1.1(a) we see that for any simple reflection $s$ in $W_0$ and $x,y$ in X, we have
$$T_s\theta_xC=\frac{(\theta_{s(x)}-\theta_x\theta_{\a_s})+v^2(\theta_x\theta_{\a_s}-\theta_{s(x)}\theta_{\a_s})}{\theta_{\a_s}-1}C,$$
$$\theta_y\theta_xC=\theta_{yx}C.$$

Thus the natural $H$-module structure on $HC$ leads  an $H$-module structure on $\CA[X]$ as follows (see [L5, Lemma 4.7]):
$$T_s\circ x=\frac{(s(x)-x\a_s)+v^2(x\a_s-s(x)\a_s)}{\a_s-1},$$
$$\theta_y\circ x=yx,$$
where $s\in W_0$ is a simple reflection and $x,y\in X$.

\def\st{\stackrel}
\def\sc{\scriptstyle}
\def\CG{\mathcal G}

\vskip.3cm \noindent{\bf 1.3.} Let $R^+$ be the set of positive roots, $R^-=-R^+$ and $\Delta$ be the set of simple roots in $R$. For $\a\in \Delta$, let $x_\a$ be the corresponding fundamental weight. For $w\in W_0$, let
$$e_w=w(\prod_{\st {\alpha\in\Delta}{ w(\alpha)\in R^-}}x_\alpha)\in X.$$
According to [St] we have

\vskip.2cm

(a) $\CA[X]$ is a free $\CA[X]^{W_0}$-module with a basis $e_w,\ w\in W_0$.

\vskip.2cm

Note that $\CA[X]^{W_0}$ is naturally isomorphic to $R_{G\times\BC^*}=R_G\otimes_{\BZ}\CA$, where $R_{\CG}$ and $R_G$ are the representation ring of ${\CG}=G\times \BC^*$ and $G$ respectively, and  we identify $R_{\BC^*}$ with $\CA$ by regarding $v$ as the identity representation $\BC^*\to\BC^*$. We also use $v$ for the element in $R_{\CG}$ defined by the natural projection $G\times\BC^*\to\BC^*.$

According to Bernstein (see [L1]), we know that

\vskip.2cm

(b) The $\theta(\CA[X]^{W_0})$ is the center $Z(H)$ of $H$ and the center $Z(H)$ of $H$ is isomorphic to the representation ring of $R_{\CG}$.

\vskip.2cm

Therefore we have

\vskip.2cm

(c) The $H$-module $HC$ is a free $Z(H)$-module with a basis $\theta_{e_w}C,\ w\in W_0$.

\vskip.3cm \noindent{\bf 1.4.} For each $w\in W_0$, set $d_w^{-1}=(\prod_{\st {\alpha\in\Delta}{ w(\alpha)\in R^-}}x_\alpha)w^{-1}\in W$. According to [Sh1, Sh2], the lowest two-sided cell of $W$ $c_0$ can be described as
$$c_0=\{ d_ww_0xd_u^{-1}\,|\, w,u\in W, \ x\in X^+\}.$$
For $x\in X^+$, let $V(x)$ be an irreducible rational $G$-module of highest weight $x$. For $x,y,z\in X^+$, let $m_{x,y,z}$ be the multiplicity of $V(z)$ appearing in $V(x)\otimes V(y)$. The based ring $J_0$ of $c_0$ defined by Lusztig in [L3] is a free $\BZ$-module with a basis $t_w,\ w\in c_0$ and the structure constants are given by (see [X1, Theorem 1.10])
$$t_{d_ww_0xd_u^{-1}}t_{d_vw_0yd_p^{-1}}=\delta_{u,v}\sum_{z\in X^+}m_{x,y,z}t_{d_ww_0zd_p^{-1}}.$$
For $x\in X^+$, let $S_x$ be the element in $\Theta$ corresponding to $V(x)$. By abuse notation, we also use $S_x$ for the element in $R_G$ or $R_{\CG}$ corresponding to $V(x)$. Then we have

\vskip.2cm

(d) The map $S_x\to \sum_{w\in W_0}t_{d_ww_0xd_w^{-1}}$ defines an ring isomorphism from $R_G$ to the center of $J_0$.

\vskip.2cm

So we have

\vskip.2cm

(e) The map $t_{d_ww_0xd_u^{-1}}\to (S_x)_{w,u}$ defines an isomorphism $\varphi_1$ of $R_G$-algebra from $J_0$ to the $W_0\times W_0$ matrix ring $M_{W_0}(R_G)$ over $R_G$, where $(S_x)_{w,u}$ stands for the $W_0\times W_0$ matrix  with entry $S_x$ at the position $(w,u)$ and with 0 entries at other positions. See [X1, Theorem 1.10].

\vskip.2cm

Note that we have $S_xC_{d_ww_0}=C_{d_ww_0x},$ see [X1, Theorem 2.9] and [L1, 8.6, 6.12]. It is known that the elements $C_{d_ww_0x},\ w\in W_0,\ x\in X^+$, form an $\CA$-basis of $HC$ (see [L2]). Hence we have

\vskip.2cm

(f) The elements $C_{d_ww_0},\ w\in W_0$, form a $Z(H)$-basis of $HC$.

\vskip.2cm

Now we have two $Z(H)$-bases for $HC$, one is $\{\theta_{e_w}C\,|\, w\in W_0\}$, the other is $\{C_{d_ww_0}\,|\, w\in W_0\}$.  For each $u\in W$, we then have
 $$C_{d_uw_0}C=\sum_{w\in W_0}a_{w,u}\theta_{e_w},\quad a_{w,u}\in Z(H)=R_{\CG}.$$ Let $A=(a_{w,u})$ be the $W_0\times W_0$ matrix. Then $A$ is the transfer matrix from the basis  $\{\theta_{e_w}C\,|\, w\in W_0\}$ to the basis $\{C_{d_ww_0}\,|\, w\in W_0\}$. That is
$$(C_{d_ww_0})_{w\in W_0}=(\theta_{e_w}C)_{w\in W_0}A,$$ where $(\cdot)_{w\in W_0}$ are row vectors.

\vskip.3cm \noindent{\bf 1.5.} For $w,u\in W$, write $C_wC_u=\sum_{v\in W}h_{w,u,v}C_v,\ h_{w,u,v}\in\CA$. Let $\mathcal D$ be the set consisting of all $d_ww_0d_w^{-1},\ w\in W_0$. Lusztig showed that the map $$\vp:H\to J_0\otimes \CA,\quad C_w\to \sum_{{\st {d\in\mathcal D}{ u\in W}}}h_{w,d,u}t_u,$$ is an $\CA$-algebra homomorphism (see [L3]). Note that the restriction of $\vp$ to the center $Z(H)$ gives an isomorphism from $Z(H)$ to the center of $J_0\otimes \CA$. If we identify $Z(H)$ with the representation ring $R_{\CG}$ of $\CG=R\times\BC^*$, then $\vp$ is in fact an $R_{\CG}$-algebra homomorphism.

\def\vp{\varphi}

Let $$\vp=(\varphi_1\otimes{\text{id}_{\CA}})\varphi_0:H\to J_0\otimes\CA\to M_{W_0}(R_{G})\otimes\CA=M_{W_0}(R_{\CG}),$$ see 1.4 (e) for the definition of $\varphi_1$. Then $\vp$ is a homomorphism of $R_{\CG}$-algebra.

\section{Equivariant K-theory Construction}

Chriss and Ginzburg constructed a homomorphism of $R_{\CG}$-algebra from $H$ to the matrix ring $M_{W_0}(R_{\CG})$, see [CG, Corollary 5.4.34, 7.6.8]. In this section we recall this construction. We shall follow the approach in [L5] for explicit calculation.

\vskip.3cm\noindent{\bf 2.1.} Let $\mathfrak{g}$ be the Lie algebra
of $G$, $\mathcal{N}$ the nilpotent cone of $\Gg$ and $\CB$  the
variety of all Borel subalgebras of $\Gg$. The Steinberg variety $Z$
is the subvariety of $\CN\times \CB\times\CB$ consisting of all
triples $(n,\gb,\gb'), \ n\in \gb\cap \gb'\cap\CN, \ \gb,\gb'\in
\CB$. Let $\Lambda=\{(n,\gb)\ | \ n\in\CN\cap \gb, \gb\in\CB\}$ be
the cotangent bundle of $\CB$. Clearly $Z$ can be regarded as a
subvariety of $\Lambda\times\Lambda$ by the imbedding $Z\to
\Lambda\times\Lambda$, $(n,\gb,\gb')\to (n,\gb,n,\gb').$ Define a
$\CG=G\times C^*$-action on $\Lambda$ by $(g,z):\ (n,\gb)\to
(z^{-2}\text{ad}(g)n,\text{ad}(g)\gb)$. Let $G\times C^*$ acts on
$\Lambda\times\Lambda$ diagonally, then  $Z$ is a $G\times
C^*$-stable subvariety of $\Lambda\times\Lambda$. Let
$K_{G\times\mathbb{C}^*}(Z)=K_{G\times\mathbb{C}^*}(\Lambda\times\Lambda;Z)$
be the Grothendieck group of the category of
$G\times\mathbb{C}^*$-equivariant coherent sheaves on
$\Lambda\times\Lambda$ with support in $Z$.

Let $$\Lambda^{aab}=\{(n,\gb,n',\gb',n'',\gb'')\in\Lambda^3\,|\,n=n'\},$$
$$\Lambda^{abb}=\{(n,\gb,n',\gb',n'',\gb'')\in\Lambda^3\,|\,n'=n''\},$$
$$\Lambda^{aaa}=\{(n,\gb,n',\gb',n'',\gb'')\in\Lambda^3\,|\,n=n'=n''\}.$$
Let $\CG$ act on $\Lambda^3$ diagonally, then $\Lambda^{aab},\  \Lambda^{abb}$ and $\Lambda^{aaa}$ are $\CG$-stable subvarieties of $\Lambda^3$.
 Define $\pi_{12}:\Lambda^{aab}\to Z,$ $\pi_{23}:\Lambda^{abb}\to Z,$ $\pi_{13}:\Lambda^{aaa}\to Z$ as follows,
 $$\pi_{12}(n,\gb,n,\gb',n,\gb'') \to (n,\gb,\gb') .$$
$$\pi_{23}(n,\gb,n',\gb',n'',\gb'') \to (n',\gb',\gb'') ,$$
$$\pi_{13}(n,\gb,n,\gb',n,\gb'') \to (n,\gb,\gb'') .$$
Note that $\pi_{12},\ \pi_{23}$ are smooth and $\pi_{13}$ is proper. Following [L5, 7.9], one defines the
convolution product
$$*:
K_{G\times\mathbb{C}^*}(Z)\times K_{G\times\mathbb{C}^*}(Z)\to
K_{G\times\mathbb{C}^*}(Z),$$
$$\mathscr{F}*\mathscr{G}=(\pi_{13})_*(\pi_{12}^*\mathscr{F}\otimes_{\mathcal
{O}_{\Lambda^3}}^Lp_{23}^*\mathscr{G}),$$
where $\mathcal {O}_{\Lambda^3}$ is the
structure sheaf of ${\Lambda^3}$. This
endows with $K_{G\times\mathbb{C}^*}(Z)$ an associative algebra
structure over the representation ring $R_{G\times \BC^*}$ of
$G\times \BC^*$. Recall that we regard the indeterminate $v$ as the
representation $G\times\BC^*\to \BC^*,\ (g,z)\to z$. Then
$R_{G\times \BC^*}$ is identified with $ R_G\otimes_{\BZ}\CA$. In
particular, $K_{G\times\mathbb{C}^*}(Z)$ is an $\CA$-algebra.
Moreover, as $\CA$-algebras, $K^{G\times\mathbb{C}^*}(Z)$ is
isomorphic to the Hecke algebra $H$, see [G1, G2, KL2] or [CG, L5]. We
shall identify $K_{G\times\mathbb{C}^*}(Z)$ with $H$.

\def\kcg{K_{\CG}}
\def\kcgz{\kcg(Z)}

\vskip.3cm \noindent{\bf 2.2.} We shall simply write $\CG$ for $G\times\BC^*$. Consider the $\CG$-equivariant injections
$$j: \CB\times\CB\to Z,\quad (\gb,\gb')\to (0,\gb,\gb'),$$
and $$k: Z\to\Lambda\times\CB,\quad (n,\gb,\gb')\to (n,\gb,\gb').$$ We then have two $K_{\CG}(Z)$-linear maps:
$$j_*: K_{\CG}(\CB\times\CB)\to K_{\CG}(Z),$$
$$k_*:\kcg(Z)\to\kcg(\Lambda\times\CB).$$
Here the $\kcgz$-module structure on $\kcg(\CB\times\CB)$ and on $\kcg(\Lambda\times\CB)$ are defined as follows.

Let
$$r_{12}:Z\times\CB\to Z,\quad (n,\gb,\gb',\gb'')\to (n,\gb,\gb'),$$
$$r_{23}:\Lambda\times\CB\times\CB\to \CB\times\CB,\quad (n,\gb,\gb',\gb'')\to (\gb',\gb''),$$
$$r_{13}:\CB\times\CB\times\CB\to \CB\times\CB,\quad (\gb,\gb',\gb'')\to (\gb,\gb''),$$
$$q_{23}:\Lambda^2\times\CB\to \Lambda\times\CB,\quad (n,\gb,n',\gb',\gb'')\to (n',\gb',\gb''),$$
$$q_{13}:Z\times\CB\to \Lambda\times\CB,\quad (n,\gb,n,\gb',\gb'')\to (n,\gb,\gb'').$$

Define the $R_{\CG}$-bilinear pairings
$$\star: \kcgz\times \kcg(\CB\times\CB)\to\kcg(\CB\times\CB),$$
$$\mathscr{F}\star\mathscr{G}=r_{13_*}(r_{12}^*\mathscr{F}\otimes^L_{\Lambda^3}r_{23}^*\mathscr{G})\in\kcg(\CB\times\CB);$$
and
$$\circ: \kcgz\times \kcg(\Lambda\times\CB)\to\kcg(\Lambda\times\CB),$$
$$\mathscr{F}\circ\mathscr{G}=q_{13_*}(r_{12}^*\mathscr{F}\otimes^L_{\Lambda^2\otimes\CB}q_{23}^*\mathscr{G})\in\kcg(\Lambda\times\CB);$$
respectively. These $R_{\CG}$-pairings give $\kcgz$-module structure on $\kcg(\CB\times\CB)$ and $\kcg(\Lambda\times\CB)$ respectively.

\vskip.3cm \noindent{\bf 2.3.}  According to [KL2,1.6] and its proof, we have
\vskip.2cm
(a) The external tensor product in $K_G$-theory
$$\boxtimes: K_G(\CB)\otimes_{R_G}K_G(\CB)\to K_G(\CB\times\CB)$$
is an isomorphism of $R_G$-module. For $x\in X$, we shall also denote the corresponding line bundle on $\CB$ by $x$.
\vskip.2cm
Note that $K_G(\CB)=R_B=\BZ[X].$ Define the pairing $(,):\BZ[X]\times\BZ[X]\to R_G$ by
$$(x,y)=\sum_{w\in W_0}(-1)^{l(w)}w(xy\rho)/\sum_{w\in W_0}(-1)^{l(w)}w(\rho),$$
where $\rho$ is the product of all fundamental weights.
\vskip.2cm
(b) There exist elements $e'_w\in\BZ[X],\ w\in W_0$, such that $(e_w,e'_u)=\delta_{w,u}$. (See the proof of [KL2,1.6].)
\vskip.2cm
 For $1\le i<j\le
3$, let $p_{ij}$ be the projection from
$\CB^3=\CB\times\CB\times\CB$ to its $(i,j)$-factor. Define the convolution $\cdot$ on $K_G(\CB\times \CB)$ as follows:
$$V\cdot V'=(p_13)_*(p_{12}^*V\otimes^L_{\CB^3}p_{23}^*V').$$
Let $w,u,t,v$ be elements in $W_0$ and $\xi,\eta$ be elements in $R_G$. Then we have
$$(e_w\xi\boxtimes e'_u)\cdot(e_{t}\boxtimes e'_v)=(e_t,e'_{u})(e_w\xi\boxtimes\eta e'_v),$$
see [X2, 5.16.1] or [CG, Lemma 5.2.28]. Therefore the map
$$ e_w\xi\boxtimes e'_u\to (\xi)_{w,u},\quad w,u\in W_0,\ \xi\in R_G,$$
defines an isomorphism of $R_G$-algebra
$$\psi_1: K_{G}(\CB\times\CB)\to M_{W_0}(R_G),$$
where $(\xi)_{w,u}$ stands for the matrix whose entries are $\xi$ at position $(w,u)$ and  zero otherwise.

\vskip.3cm \noindent{\bf 2.4.} The projection $pr:\Lambda\times\CB\to\CB\times\CB,$ $(n,\gb,\gb')\to(\gb,\gb')$, is $\CG$-equivariant. The Thom isomorphism says that $pr^*:\kcg(\CB\times\CB)\to\kcg(\Lambda\times\CB)$ is isomorphism of $R_{\CG}$-module. Composing its inverse with $k_*:\kcgz\to\kcg(\Lambda\times\CB)$ we get an $R_{\CG}$-linear map
$$(pr^*)^{-1}k_*: \kcgz\to\kcg(\CB\times\CB).$$ The following result is showed in [CG, Corollary 5.4.34].

\vskip.2cm

(a) The map $(pr^*)^{-1}k_*: \kcgz\to\kcg(\CB\times\CB)$ is a homomorphism of $R_{\CG}$-algebra.

\vskip.2cm

We explain this. Through inverse image $pr_2^*:\kcg(\CB)\to \kcg(\Lambda)$ of the projection $pr_2:\Lambda\to  \CB$, we regard the trivial bundle $\BC$  on $\CB$ as an $\CG$-equivariant bundle on $\Lambda$. Let $r:\Lambda\to Z$ be the inclusion $(n,\gb)\to (n,\gb,\gb)$. Then $r_*(\BC)$ is the unit in $\kcgz$ (see [L5, 7.10].

Let $\tilde j: \Lambda\to \Lambda\times \CB$ be the inclusion $(n,\gb)\to(n,\gb,\gb)$ and let $i:\CB\to \CB\times\CB$ be the diagonal map. Then we have a cartesian  diagram
 $$\xymatrix{
  \Lambda \ar[d]_{pr_2} \ar[r]^{\tilde j}
                & \Lambda\times\CB \ar[d]^{pr}  \\
  \CB \ar[r]_{i}
                & \CB\times\CB             }$$
So we have $\tilde j_*pr_2^*(\BC)=pr^*i_*(\BC)$ (see for example [FK, Chapter I, 6.1 Theorem]). Therefore, $(pr^*)^{-1}\tilde j_*pr_2^*(\BC)=i_*(\BC)$. This means that $(pr^*)^{-1}k_*r_*(\BC)$ is  $i_*(\BC)$ in $\kcg(\CB\times\CB)$. But it is known that $i_*(\BC)=\sum_{w\in W_0}e_w\boxtimes e'_w$ is the unit in $K_{\CG}(\CB\times\CB)$ (see [KL2, 1.7]). Since $k_*$ is $\kcgz$-linear, we see easily that the map in (a) is a homomorphism of $R_{\CG}$-algebra.

Define $\psi: H=\kcgz\to M_{W_0}(R_{\CG})=M_{W_0}(R_G)\otimes \CA$ to be the homomorphism of $R_{\CG}$-algebra $(\psi_1\otimes\text{id}_{\CA})(pr^*)^{-1}k_*$, see 2.3 for definition of $\psi_1$. Recall that  we defined in subsection 1.5 the homomorphism $\vp:H\to M_{W_0}(R_{\CG})$ of $R_{\CG}$-algebra and defined the matrix $A\in M_{W_0}(R_{\CG})$ in the last part of subsection 1.4. Now we can state the main result of the paper.

\vskip.3cm \noindent{\bf Theorem 2.5.} For $h\in H$ we have $\vp(h)=A^{-1}\psi(h)A$.

Proof. Let $h\in H$. According to [L5, 7.14], we have $$(pr^*)^{-1}k_*(h)=(pr^*)^{-1}k_*(h)(pr^*)^{-1}k_*(1).$$ In subsection 2.4 we have seen that $$(pr^*)^{-1}k_*(1)=\sum_{w\in W_0}e_w\boxtimes e'_w.$$ According to [L5, Lemma 7.21], for $h\in H$, we have  $$(pr^*)^{-1}k_*(h)=h\circ \sum_{w\in W_0}e_w\boxtimes e'_w=\sum_{w\in W_0}(h\circ e_w)\boxtimes e'_w.$$

Let $h\circ e_w=\sum_{u\in W_0}h_{u,w}e_u.$ Then we have $$\psi(h)=(h_{u,w})\in M_{W_0}(R_{\CG}).$$ Note that we also have $$h\theta_{e_w}C=\sum_{u\in W_0}h_{u,w}\theta_{e_u}C.$$ Therefore the homomorphism comes from the natural module structure on $HC$ of the $R_{\CG}$-algebra $H$ with respect to the $R_\CG$-basis $\theta_{e_w}C,\ w\in W_0$.

By definition, for $r\in W$, we have $$\varphi_0(C_r)=\sum_{d\in \CD}h_{r,d,r'}t_{r'},\ h_{r,d,r'}\in \CA.$$
By [L2] and [Sh2], $h_{r,d_ww_0d_w^{-1},r'}\ne 0$ if and only if $r'=d_uw_0xd_w^{-1}$ for some $u\in W_0$ and $x\in X^+$. Using [X1, Theorem 2.9] we know that  $h_{r,d_ww_0d_w^{-1},r'}=h_{r,d_ww_0,d_uw_0x}$ for $r'=d_uw_0xd_w^{-1}$, $u,w\in W_0,$ $x\in X^+$. Moreover, for $x\in X^+$, we have $S_xC_{d_uw_0d_w^{-1}}=C_{d_uw_0xd_w^{-1}}$ and $S_xC_{d_uw_0}=C_{d_uw_0x},$ see [X1, Theorem 2.9] and [L1, 8.6, 6.12]. Therefore, for $h\in H$, if $\varphi(h)=(b_{u,w})\in M_{W_0}(R_{\CG})$, then we have $hC_{d_ww_0}=\sum_{u,w}b_{u,w}C_{d_uw_0}$. Hence, the homomorphism $\vp$ also comes from the natural module structure on $HC$ of the $R_{\CG}$-algebra $H$, but with respect to the $R_\CG$-basis $C_{d_ww_0},\ w\in W_0$.

 The theorem follows.

\vskip.3cm{\bf 2.6.} According to [L5, 10.6], we have
$$j_*(x\boxtimes y)=(-1)^\nu\theta_{x\rho}\sum_{w\in W_0}T_w\theta_{\rho y},$$ where $\nu=|R^+|$. By the proof of Proposition 10.7 in [L5], we see that
$$(pr^*)^{-1}k_*j_*: K_{\CG}(\CB\times\CB)\to\kcgz\to K_{\CG}(\Lambda\times \CB)\to K_{\CG}(\CB\times\CB)$$ is given by the following formula
$$(pr^*)^{-1}k_*j_*(x\boxtimes y)=(\prod_{\alpha\in R^+}(1-v^2\alpha)) x\boxtimes y.$$

Let $x=y=\rho^{-1}$. In view of the discussion in 1.2 and [L5, Lemma 7.21], the above formula implies that
$$\sum_{w\in W_0}T_w\theta_{\rho^{-1}} C=\prod_{\alpha\in R^+}(1-v^2\theta_\alpha)\theta_{\rho^{-1}}C.$$
This formula is proved in [L1], now has an interpretation in terms of equivariant K-theory.


\vskip3mm
 {\bf Acknowledgement:} The work was completed during the author's visit to the Institute Mittag-Leffler. The author is grateful to the Institute for hospitality. Part of the work was done during the  author's visit to the Departement de Mathematiques, Universite de Paris VII in 2012. The author is very grateful to E. Vasserot for invitation and for very helpful discussions and comments.

\bibliographystyle{unsrt}

\end{document}